\documentclass[english]{article}
\usepackage[T1]{fontenc}
\usepackage[latin9]{inputenc}
\usepackage{amsmath}
\usepackage{graphicx}
\usepackage{amssymb}
\usepackage{esint}

\makeatletter

\providecommand{\tabularnewline}{\\}

\makeatother

\usepackage{babel}

\begin{document}

\title{Exact \& Numerical Tests of Generalised Root Identities for non-integer
$\mu$}

\author{Richard Stone}
\maketitle
\begin{abstract}
We consider the generalised root identities introduced in {[}1{]}
for simple functions, and also for $\Gamma(z+1)$ and $\zeta(s)$.
Where these were almost exclusively examined just for $\mu\in\mathbb{\mathbb{Z}}$
in {[}1{]}, here we focus exclusively on the case of $\mu\mathbb{\mathbb{\notin Z}}$.
For the simplest function $f(z)=z$, and hence for arbitrary polynomials,
we show that they are satisfied for arbitrary $\mu\in\mathbb{R}$
(and hence for arbitrary $\mu\in\mathbb{\mathbb{C}}$ by analytic
continuation). Using this, we then develop an asymptotic formula for
the derivative side of the root identities for $\Gamma(z+1)$ at arbitrary
$\mu\in\mathbb{R}$, from which we are able to demonstrate numerically
that $\Gamma(z+1)$ also satisfies the generalised root identities
for arbitrary $\mu$, not just integer values as verified in {[}1{]}.
Finally, using the working in {[}1{]} we examine the generalised root
identites for $\zeta$ also for non-integer values of $\mu${\normalsize .}
Having shown in {[}1{]} that $\zeta$ satisfies these identities exactly
for $\mu\in\mathbb{Z_{\geq\textrm{2}}}$ (and also for $\mu=1$ after
removal of an obstruction), in this paper we present strong numerical
evidence first that $\zeta$ satisfies them for arbitrary $\mu>1$
where the root side is classically convergent, and then that this
continues to be true also for $-1<\mu<1$ where Cesaro divergences
must be removed and Cesaro averaging of the residual partial-sum functions
is required (when $\mu<0$). A careful examination of the calculations
in a neighbourhood of $\mu=0$ also sheds light on the appearance
of the 2d ln-divergence that was handled heuristically in {[}1{]}
and why the assignment of 2d Cesaro limit $0$ to this in {[}1{]}
is justified. The numerical calculations for $\mu>0$ are bundled
in portable R-code which can be readily used to further verify these
claims; the code for the case $-1<\mu<0$, including the Cesaro averaging
required when $\mu<0$, is in VBA in an XL spreadsheet which can also
be used to check these claims and conduct further tests. Both the
R-scripts and XL spreadsheet are made available with this paper, along
with supporting files.
\end{abstract}

\section{Introduction}

In {[}1{]} we introduced the generalised root identities, namely:

\begin{equation}
\frac{-1}{\Gamma(\mu)}\left(\frac{\textrm{d}}{\textrm{d}z}\right)^{\mu}(\ln(f(z)))|_{z=z_{0}}=\textrm{e}^{i\pi\mu}\sum_{\{z_{0}-roots\, r_{i}\}}\frac{M_{i}}{(z_{0}-r_{i})^{\mu}}\label{eq:GenRtId1}\end{equation}
The LHS in (\ref{eq:GenRtId1}) is the derivative side of the root
identities for $f$ at $z_{0}$ and $\mu$; the RHS is the root side.
Writing $d(s_{0},\mu)$ for the expression on the derivative side
and $r(s_{0},\mu)$ for the root side, saying that a function $f$
satisfies the generalised root identities means that $d(s_{0},\mu)$
and $r(s_{0},\mu)$ are identically equal as functions of two complex
variables (we may omit reference to the underlying function $f$ in
the notation for $d$ and $r$ unless the context requires it). 

In general, within an equivalence class of functions (equivalent up
to multiplication by a nowhere-zero entire function), we expect there
to be a unique representative which satisfies these root identities.
For many functions of great interest, however, we saw in {[}1{]} that
they themselves satisfy the generalised root identities, at least
for many values of $s_{0}$ and $\mu$, without any need to go exploring
within their equivalence classes. For example in {[}1{]} we showed
that:

(a) the function $f(z)=\cos(\frac{\pi z}{2})$ satisfies the generalised
root identities at $z_{0}=0$ for arbitrary $\mu\in\mathbb{\mathbb{C}}$,
with the root identities in this case being precisely equivalent to
the famous functional equation for $\zeta$ in $\mu$.

(b) the function $\Gamma(z+1)$ satisfies the generalised root identities
for arbitrary $z_{0}\in\mathbb{C}$ and $\mu\in\mathbb{Z}$, with
a renormalisation being required when $\mu=1$, and with a distributional
interpretation of the derivative side and a Cesaro evaluation of the
root side required when $\mu\in\mathbb{Z_{\leq\textrm{0}}}$.

(c) the function $\zeta(s)$ satisfies the generalised root identities
for arbitrary $s_{0}\in\mathbb{C}$ and $\mu\in\mathbb{Z_{\geq\textrm{2}}}$,
while $\pi^{-\frac{s}{2}}\zeta(s)$ satisfies the identities when
$\mu=1$ after the removal of an {}``obstruction'' in this instance
(the resulting $\mu=1$ identity being equivalent to the Hadamard
product formula for $\zeta$). Additionally we showed that the derivative
side of (\ref{eq:GenRtId1}) should be identically zero as a function
of $s_{0}$ when $\mu\in\mathbb{Z_{\leq\textrm{0}}}$ and we conducted
three Cesaro calculations of the root side for such cases, namely
for $\mu=0$, $\mu=-1$ and $\mu=-2$. For $\mu=0$ and $\mu=-1$
we verified the claimed result, at least modulo an estimate on the
argument of the zeta function, $S(T)$, in the case $\mu=-1$ which
is known to hold if the Riemann hypothesis (RH) is true. For $\mu=-2$,
however, we claim that the resulting Cesaro evaluation of the root
side suggests that the Riemann hypothesis must be false.

As discussed in {[}1{]}, however, the fact that these computations
were conducted for isolated values of $\mu$, rather than as a systematic
extension in open regions of the $\mu$-plane, is a deficiency of
the argumentation thus far (in {[}1{]}, the case of $f(z)=\cos(\frac{\pi z}{2})$
when $z_{0}=0$ is in fact the only instance in the whole paper where
the case of arbitrary $\mu\notin\mathbb{Z}$ is considered for root
identity calculations).

In this paper, therefore, we address this by explicitly focusing on
the root identities for arbitrary $\mu$ and $z_{0}$ in $\mathbb{R}$
(and then $\mathbb{C}$ by analytic continuation), verifying them
in this general case also for the functions considered in {[}1{]},
in particular $\Gamma$ and $\zeta$.

As in {[}1{]}, however, we start with simple examples before proceeding
to the more complex ones. Specifically, in section 2.1, we first consider
the simplest possible example of $f(z)=z$ and show that it satisfies
the generalised root identities (\ref{eq:GenRtId1}) for arbitrary
$(z_{0},\mu)\in\mathbb{C^{\textrm{2}}}$. In this case, of course,
the root side, $r(s_{0},\mu)$, is trivial to compute, but the derivative
side for arbitrary $\mu$ is already interesting, requiring the usual
distributional calculations when $\mu\in\mathbb{Z_{\leq\textrm{0}}}$
and the use of generalised Fresnel integrals when $\mu\in\mathbb{R\setminus\mathbb{Z}}$.

We then use this result to extend to the case of arbitrary polynomials
in section 2.2, and to the case of $\Gamma(z+1)$ in section 3.

Finally we consider the case of $\zeta$ in section 4. By considering
the case of arbitrary $\mu$ and $s_{0}$ we shall be able to verify
numerically that $\zeta$ satisfies the generalised root identities
(\ref{eq:GenRtId1}) when $\mu$ and $s_{0}$ are in half-planar regions
of $\mathbb{C}$, and also discuss how these numerical computations
may be both extended to all of $(s_{0},\mu)\in\mathbb{C^{\textrm{2}}}$
and used to cast light on certain other theoretical gaps identified
in {[}1{]}, in particular the legitimacy of the claim there that we
may take $\underset{z,\tilde{z}\rightarrow\infty}{Clim}\ln\left(\frac{z}{\tilde{z}}\right)=0$
on the root side for $\zeta$ for $\mu=0$, $\mu=-1$ and $\mu=-2$.
These results in turn, we believe, give strong numerical evidence
in support of the idea that $\zeta$ does satisfy the generalised
root identities at arbitrary $\mu$ (and in particular at $\mu=0,-1$,
and $-2$) and hence that the claims in {[}1{]} are valid.

\section{The Root Identities for $f(z)=z$ and for Polynomials}

\subsection{The Case of $f(z)=z$}

In this case $f$ has a single root of multiplicity $1$ at $z=0$
and so the root-side of the generalised root identities is given by

\begin{equation}
r(z_{0},\mu)=\frac{\textrm{e}^{i\pi\mu}}{z_{0}^{\mu}}\label{eq:RtId_id_fn_RootSide}\end{equation}

We now consider the derivative side in several steps.

(i) The case of $\mu\in\mathbb{Z_{>\textrm{0}}}$: When $\mu=1$ the
derivative side is 

\begin{equation}
d(z_{0},1)=\frac{-1}{\Gamma(1)}\frac{\textrm{d}}{\textrm{d}z}(\ln z)|_{z=z_{0}}=\frac{-1}{z_{0}}=r(z_{0},1)\label{eq:RtId_id_fn_DerivSide1}\end{equation}
and by direct differentiation it is then trivial to see that we likewise
have

\begin{equation}
d(z_{0},\mu)=\frac{\textrm{e}^{i\pi\mu}}{z_{0}^{\mu}}=r(z_{0},\mu)\quad\quad\textrm{f\textrm{or\;\ all}}\;\mu\in\mathbb{Z_{>\textrm{0}}}\label{eq:RtId_id_fn_DerivSide2}\end{equation}
so that certainly $f(z)=z$ satisfies the generalised root identities
for $\mu\in\mathbb{Z_{>\textrm{0}}}$.%
\footnote{Note in passing that these results can be derived via integration
from the Fourier definition of the derivative side in {[}1{]} using
Cesaro methods, but there is obviously no need when $\mu\in\mathbb{Z_{>\textrm{0}}}$%
}

(ii) The case of $\mu\in\mathbb{Z_{\leq\textrm{0}}}$: Recalling the
definition from {[}1{]} of the derivative side, that\begin{equation}
d_{f}(z_{0},\mu)=-\frac{1}{\Gamma(\mu)}\frac{1}{2\pi}\int_{-\infty}^{\infty}(i\xi)^{\mu}\mathcal{F}[\ln f](\xi)\textrm{e}^{iz_{0}\xi}\textrm{d}\xi\label{eq:DerivSideDef}\end{equation}
and that $\mathcal{F}[\ln z](\xi)=-2\pi$$\frac{\tilde{H}_{0}(\xi)}{\xi}$
(where $\tilde{H}_{0}(\xi)$ is the odd Heaviside function) we can
calculate the derivative side via the sort of distributional calculations
used in {[}1{]}. For $\mu=0$ we have

\begin{equation}
d(z_{0},0)=\frac{1}{\Gamma(0)}\int_{-\infty}^{\infty}\frac{\tilde{H}_{0}(\xi)}{\xi}\,\textrm{e}^{iz_{0}\xi}\textrm{d}\xi=\delta_{0}(\xi)[\textrm{e}^{iz_{0}\xi}]=1=r(z_{0},0)\label{eq:RtId_id_fn_DerivSide3}\end{equation}
while for $\mu=-1$ we have

\begin{equation}
d(z_{0},-1)=\frac{-i}{\Gamma(-1)}\int_{-\infty}^{\infty}\frac{\tilde{H}_{0}(\xi)}{\xi^{2}}\,\textrm{e}^{iz_{0}\xi}\textrm{d}\xi=i\cdot\delta_{0}^{\prime}(\xi)[\textrm{e}^{iz_{0}\xi}]=-z_{0}=r(z_{0},-1)\label{eq:RtId_id_fn_DerivSide4}\end{equation}
and similarly for $\mu=-2,-3,\ldots$ so that $f(z)=z$ also satisfies
the generalised root identities for $\mu\in\mathbb{Z_{\leq\textrm{0}}}$.

(iii) The case of $\mu\in(0,1)$: Assume initially that $z_{0}$ is
also real and positive. Then the derivative side, $d(z_{0},\mu)$
is given by

\begin{eqnarray}
d(z_{0},\mu) & = & \frac{1}{\Gamma(\mu)}\int_{-\infty}^{\infty}(i\xi)^{\mu}\frac{\tilde{H}_{0}(\xi)}{\xi}\textrm{e}^{iz_{0}\xi}\textrm{d}\xi\label{eq:RtId_id_fn_DerivSide5a}\\
\nonumber \\ & = & \frac{\textrm{e}^{i\frac{\pi}{2}\mu}}{\Gamma(\mu)}\int_{-\infty}^{\infty}\tilde{H}_{0}(\xi)\xi^{\mu-1}(\cos(z_{0}\xi)+i\sin(z_{0}\xi))\textrm{d}\xi\nonumber \\
\nonumber \\ & = & \frac{\textrm{e}^{i\frac{\pi}{2}\mu}}{2\Gamma(\mu)}\left\{ \begin{array}{cc}
(1-\textrm{e}^{i\pi(\mu-1)})\int_{0}^{\infty}\xi^{\mu-1}\cos(z_{0}\xi)\textrm{d}\xi\\
\\+i\cdot(1+\textrm{e}^{i\pi(\mu-1)})\int_{0}^{\infty}\xi^{\mu-1}\sin(z_{0}\xi)\textrm{d}\xi\end{array}\right\} \label{eq:RtId_id_fn_DerivSide5}\end{eqnarray}
on recalling the definition of $\tilde{H}_{0}(\xi)$. Changing variables
in these integrals by letting $v=(z_{0}\xi)^{\mu}$ and noting that
$\textrm{d}v=\mu z_{0}^{\mu}\xi^{\mu-1}\textrm{d}\xi$, it follows
that 

\begin{equation}
d(z_{0},\mu)=\frac{\textrm{e}^{i\frac{\pi}{2}\mu}}{2\Gamma(\mu+1)}\left\{ \begin{array}{cc}
(1+\textrm{e}^{i\pi\mu})\int_{0}^{\infty}\cos(v^{\frac{1}{\mu}})\textrm{d}v\\
\\+i\cdot(1-\textrm{e}^{i\pi\mu})\int_{0}^{\infty}\sin(v^{\frac{1}{\mu}})\textrm{d}v\end{array}\right\} \cdot z_{0}^{-\mu}\label{eq:RtId_id_fn_DerivSide6}\end{equation}

But recall the generalised Fresnel integral identity

\begin{equation}
\int_{0}^{\infty}\sin(x^{a})\textrm{d}x=\frac{\Gamma(\frac{1}{a})\sin(\frac{\pi}{2a})}{a}\label{eq:GeneralisedFresnel}\end{equation}
defined initially for $a>1$ and then extended by analytic continuation.%
\footnote{When $a=2$ this reduces to the classical Fresnel integral from optics
$\int_{0}^{\infty}\sin(x^{2})\textrm{d}x=\surd\frac{\pi}{8}$%
}

From this it readily follows also that

\begin{equation}
\int_{0}^{\infty}\cos(x^{a})\textrm{d}x=\frac{\Gamma(\frac{1}{a})\cos(\frac{\pi}{2a})}{a}\label{eq:GeneralisedFresnel2}\end{equation}
for $a>1$ and extended by analytic continuation.

Combining (\ref{eq:GeneralisedFresnel}) and (\ref{eq:GeneralisedFresnel2})
in (\ref{eq:RtId_id_fn_DerivSide6}) it follows that we have

\begin{eqnarray}
d(z_{0},\mu) & = & \frac{\textrm{e}^{i\frac{\pi}{2}\mu}\mu\Gamma(\mu)}{2\Gamma(\mu+1)}\left\{ (1+\textrm{e}^{i\pi\mu})\cos(\frac{\pi\mu}{2})+i\cdot(1-\textrm{e}^{i\pi\mu})\sin(\frac{\pi\mu}{2})\right\} \cdot z_{0}^{-\mu}\nonumber \\
\nonumber \\ & = & \frac{\textrm{e}^{i\frac{\pi}{2}\mu}}{4}\left\{ (1+\textrm{e}^{i\pi\mu})(\textrm{e}^{i\frac{\pi}{2}\mu}+\textrm{e}^{-i\frac{\pi}{2}\mu})+(1-\textrm{e}^{i\pi\mu})(\textrm{e}^{i\frac{\pi}{2}\mu}-\textrm{e}^{-i\frac{\pi}{2}\mu})\right\} \cdot z_{0}^{-\mu}\nonumber \\
\nonumber \\ & = & \frac{\textrm{e}^{i\pi\mu}}{z_{0}^{\mu}}=r(z_{0},\mu)\label{eq:RtId_id_fn_DerivSide7}\end{eqnarray}

Thus we see that $f(z)=z$ also satisfies the generalised root identities
for arbitrary $\mu\in(0,1)$ and $z_{0}\in\mathbb{R_{>\textrm{0}}}$. 

(iv) Extension to arbitrary $\mu,z_{0}\in\mathbb{C}$: Keeping $z_{0}\in\mathbb{R_{>\textrm{0}}}$
initially, we can then extend from $\mu\in(0,1)$ to all of $\mu\in\mathbb{R\setminus\mathbb{Z}}$
by a combination of differentiation under the integral w.r.t $z_{0}$
and homogeneity arguments (which guarantee there are no constants
of integration).

We thus have now verified that $f(z)=z$ satisfies the generalised
root identities for arbitrary $\mu\in\mathbb{R}$ and $z_{0}\in\mathbb{R_{>\textrm{0}}}$.
By unique analytic continuation they therefore continue to hold first
for arbitrary $\mu\in\mathbb{\mathbb{C}}$ and $z_{0}\in\mathbb{R_{>\textrm{0}}}$,
and finally also for arbitrary $\mu\in\mathbb{\mathbb{C}}$ and $z_{0}\in\mathbb{C}$,
except for a branch cut from $z_{0}=-\infty$ to $z_{0}=0$ when $\mu\notin\mathbb{Z}$.

\subsection{The Case of Polynomials}

Since $f(z)=z$ satisfies the generalised root identities for arbitrary
$\mu,z_{0}\in\mathbb{\mathbb{C}}$ it follows readily from (\ref{eq:DerivSideDef})
and elementary properties of the Fourier transform (in particular
that $\mathcal{F}[f(x+a)](\xi)=\mathcal{F}[f](\xi)\textrm{e}^{ia\xi}$)
that the generalised root identities are also satisfied for arbitrary
$z_{0},\mu\in\mathbb{\mathbb{C}}$ by $f(z)=z-a$.

But then it also follows immediately from the fact that $\ln((z-r_{1})(z-r_{2}))=\ln(z-r_{1})+\ln(z-r_{2})$
and the obvious additivity of the root side that the generalised root
identities must also be satisfied for arbitrary $z_{0},\mu\in\mathbb{\mathbb{C}}$
by any polynomial $p(z)$. For any given $p(z)$, the structure of
the branch cuts in $z_{0}$ when $\mu\notin\mathbb{Z}$ (corresponding
to the branch cut on $z_{0}\in(-\infty,0)$ for $f(z)=z$ replicated
for each factor) is slightly more complicated, but there are no branch
cut complications for generic $z_{0}$, and in particular for all
$z_{0}$ in a half-plane to the right of the root with the most positive
real part, which is all we require in certain subsequent arguments.

\section{The Root Identities for $\Gamma(z+1)$}

In this case $\Gamma(z+1)$ has poles of multiplicity $M_{i}=-1$
at $-1,-2,-3,\ldots$ and so the root side of the generalised root
identities is

\begin{equation}
r_{\Gamma}(z_{0},\mu)=-\textrm{e}^{i\pi\mu}\sum_{\{z_{0}+j\}}\frac{1}{(z_{0}+j)^{\mu}}\label{eq:GammaRootSide}\end{equation}
which is classically convergent for $Re(\mu)>1$ and Cesaro convergent
for general $\mu$ (renormalisation being required as per {[}1{]}
when $\mu=1$). For example, using the Euler-McLaurin sum formula
and the same sort of reasoning as used in {[}1{]}, when $0<\mu<1$
we have $r_{\Gamma}(z_{0},\mu)$ given by

\begin{equation}
r_{\Gamma}(z_{0},\mu)=-\textrm{e}^{i\pi\mu}\underset{k\rightarrow\infty}{lim}\left\{ \sum_{j=1}^{k}\frac{1}{(z_{0}+j)^{\mu}}-\frac{(z_{0}+k)^{1-\mu}}{1-\mu}\right\} \label{eq:GammaRootSide1}\end{equation}
while for $-1<\mu<0$ we have

\begin{equation}
r_{\Gamma}(z_{0},\mu)=-\textrm{e}^{i\pi\mu}\underset{k\rightarrow\infty}{lim}\left\{ \sum_{j=1}^{k}\frac{1}{(z_{0}+j)^{\mu}}-\frac{(z_{0}+k)^{1-\mu}}{1-\mu}-\frac{1}{2}(z_{0}+k)^{-\mu}\right\} \label{eq:GammaRootSide2}\end{equation}
and so on. For $\mu\notin\mathbb{Z}$ there is a single branch cut
on $z_{0}\in(-\infty,-1)$ (with branch points at $-1,-2,-3,\ldots$).

On the derivative side, using the working from {[}1{]}, section 3.4,
example 2, we have\begin{eqnarray}
d_{\Gamma}(z_{0},\mu) & = & \frac{1}{\Gamma(\mu)}\int_{-\infty}^{\infty}\tilde{H}_{0}(\xi)\,(i\xi)^{\mu-1}\left(i\frac{\textrm{e}^{i\xi}}{\textrm{e}^{i\xi}-1}\right)\textrm{e}^{iz_{0}\xi}\textrm{d}\xi\nonumber \\
\nonumber \\ & = & -\frac{i\textrm{e}^{i\frac{\pi}{2}\mu}}{\Gamma(\mu)}\int_{-\infty}^{\infty}\tilde{H}_{0}(\xi)\,\xi^{\mu-2}\left\{ \begin{array}{cc}
1-\left[\frac{(i\xi)}{2!}+\frac{(i\xi)^{2}}{3!}+\ldots\right]+\\
\\\left[\frac{(i\xi)}{2!}+\frac{(i\xi)^{2}}{3!}+\ldots\right]^{2}-\ldots\end{array}\right\} \textrm{e}^{i(z_{0}+1)\xi}\textrm{d}\xi\nonumber \\
\label{eq:GammaDerivSide1}\end{eqnarray}

In {[}1{]} we showed already that $d_{\Gamma}(z_{0},\mu)=r_{\Gamma}(z_{0},\mu)$
when $\mu\in\mathbb{Z}$ (after renormalisation when $\mu=1$) and
thus that $\Gamma$ satisfies the generalised root identities when
$\mu$ is an integer.

Now, taking the case when $\mu\in\mathbb{R\setminus\mathbb{Z}}$,
the working in section 2.1, in particular (\ref{eq:RtId_id_fn_RootSide})
and (\ref{eq:RtId_id_fn_DerivSide5a}) has shown that

\begin{equation}
\frac{\textrm{e}^{i\frac{\pi}{2}\rho}}{\Gamma(\rho)}\int_{-\infty}^{\infty}\tilde{H}_{0}(\xi)\,\xi{}^{\rho-1}\textrm{e}^{iz_{0}\xi}\textrm{d}\xi=\frac{\textrm{e}^{i\pi\rho}}{z_{0}^{\rho}}\label{eq:GammaDerivSide2}\end{equation}
for arbitrary $\rho\in\mathbb{R\setminus\mathbb{Z}}$ and $z_{0}\notin(-\infty,0)$.

It follows in (\ref{eq:GammaDerivSide1}) that, for arbitrary $\mu\in\mathbb{R\setminus\mathbb{Z}}$
and $z_{0}\notin(-\infty,-1)$, we have, after simplification,

\begin{eqnarray}
d_{\Gamma}(z_{0},\mu) & = & \left\{ \begin{array}{cc}
\frac{1}{(\mu-1)}\textrm{e}^{i\pi(\mu-1)}\frac{1}{(z_{0}+1)^{(\mu-1)}}-\frac{1}{2}\textrm{e}^{i\pi\mu}\frac{1}{(z_{0}+1)^{\mu}}\\
\\+\frac{\mu}{12}\textrm{e}^{i\pi(\mu+1)}\frac{1}{(z_{0}+1)^{(\mu+1)}}\\
\\-\frac{\mu(\mu+1)(\mu+2)}{720}\textrm{e}^{i\pi(\mu+3)}\frac{1}{(z_{0}+1)^{(\mu+3)}}+\ldots\end{array}\right\} \label{eq:GammaDerivSide3}\\
\nonumber \\ & = & -\frac{\textrm{e}^{i\pi\mu}}{(\mu-1)}\frac{1}{(z_{0}+1)^{(\mu-1)}}\left\{ \begin{array}{cc}
1+\frac{1}{2}\frac{\mu-1}{(z_{0}+1)}\\
\\+\frac{(\mu-1)\mu B_{2}}{2!(z_{0}+1)^{2}}\\
\\+\frac{(\mu-1)\mu(\mu+1)(\mu+2)B_{4}}{4!(z_{0}+1)^{4}}+\ldots\end{array}\right\} \label{eq:GammaDerivSide3a}\end{eqnarray}
as an Euler-McLaurin-style asymptotic expansion for the derivative
side. 

Using (\ref{eq:GammaRootSide})-(\ref{eq:GammaDerivSide3a}) we may
then write code to numerically check the generalised root identities
for $\Gamma(z+1)$ when $\mu\notin\mathbb{Z}$. In particular, for
any given $\mu\notin\mathbb{Z}$ we can evaluate $d_{\Gamma}(z_{0},\mu)$
to any desired degree of accuracy for $\Re(z_{0})$ sufficiently large
using sufficient terms in (\ref{eq:GammaDerivSide3}); and verifying
$d_{\Gamma}(z_{0},\mu)=r_{\Gamma}(z_{0},\mu)$ on such a half-plane
in $z_{0}$ (sufficiently far to the right for the given $\mu$) suffices
to verify their agreement for all $z_{0}$ by analytic continuation
for the given $\mu$.

Numerical investigations implementing (\ref{eq:GammaRootSide})-(\ref{eq:GammaDerivSide3a})
have been implemented in R-code by Dana Pascovici and it has been
found that $\Gamma(z+1)$ does indeed appear to satisfy the generalised
root identities for arbitrary $z_{0},\mu\in\mathbb{\mathbb{C}}$ (with
branch cut on $z_{0}\in(-\infty,-1)$ when $\mu\notin\mathbb{Z}$).
For example, consider the following:

\includegraphics[scale=0.5]{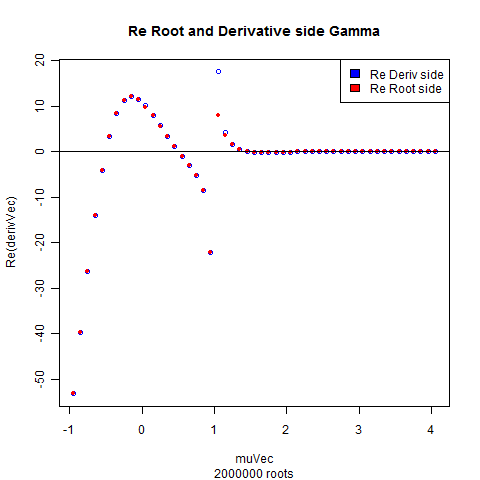}

\includegraphics[scale=0.5]{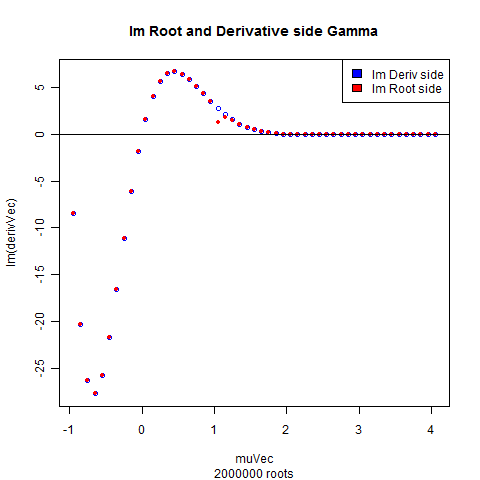}

These pictures show agreement between real and imaginary parts of
$d_{\Gamma}(z_{0},\mu)$ and $r_{\Gamma}(z_{0},\mu)$ for a range
of $\mu$-values from $-1$ to $4$ for $z_{0}=10.381$ using $2,000,000$
roots to approximate $r_{\Gamma}(z_{0},\mu)$ and truncating the expression
for $d_{\Gamma}(z_{0},\mu)$ at $\frac{1}{(z_{0}+1)^{(\mu+5)}}$.
The following picture also shows how, for $\mu$ small ($\mu=-0.01$),
the function $r_{\Gamma}(z_{0},-0.01)$ closely approaches the limiting
function $r_{\Gamma}(z_{0},0)$ which we know from section 3.4 in
{[}1{]} is given by $r_{\Gamma}(z_{0},0)=z_{0}+\frac{1}{2}$. 

\includegraphics[scale=0.47]{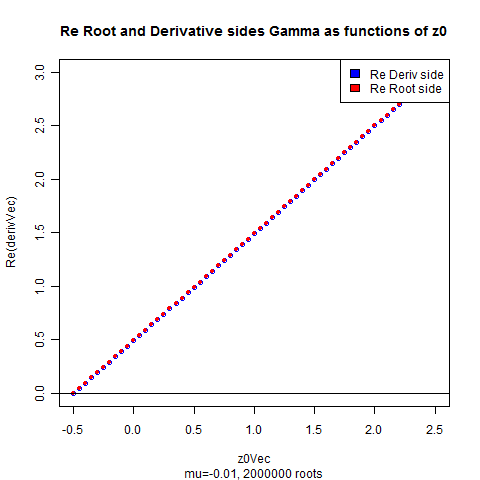}

The R-script for Gamma is made available with this paper and further
numerical tests may be performed as desired.

\section{The Root Identities for $\zeta(s)$}

In {[}1{]}, section 4.1, we derived an expression, using the Euler
product formula, for the derivative side of the generalised root identities
for $\zeta$, which applies for arbitrary $s_{0},\mu\in\mathbb{\mathbb{C}}$,
namely

\begin{equation}
d_{\zeta}(s_{0},\mu)=-\frac{\textrm{e}^{i\pi\mu}}{\Gamma(\mu)}\;\sum_{p\, prime}(\ln p)^{\mu}\left\{ \begin{array}{cc}
p^{-s_{0}}+2^{\mu-1}p^{-2s_{0}}\\
\\+3^{\mu-1}p^{-3s_{0}}+4^{\mu-1}p^{-4s_{0}}+\ldots\end{array}\right\} \label{eq:zetaDerivSideA}\end{equation}
For $\Re(s_{0})>1$ this is convergent for arbitrary $\mu\in\mathbb{C}$;
and for any given $\mu$ can then be uniquely analytically continued
to arbitrary $s_{0}\in\mathbb{\mathbb{C}}$ from this right half-plane
in $s_{0}$. In {[}1{]} we confirmed that $\zeta$ satisfies the generalised
root identities when $\mu\in\mathbb{Z_{\geq\textrm{2}}}$ (for which
the root side is convergent) and that it satisfies them for $\mu=1$
also after removal of an obstruction, with this result for $\mu=1$
corresponding to the Hadamard product formula for the related function
$\xi$. We also noted that (for $\Re(s_{0})>1$ initially and then
all $s_{0}\in\mathbb{C}$ by analytic continuation) equation (\ref{eq:zetaDerivSideA})
clearly implies that 

\begin{equation}
d_{\zeta}(s_{0},\mu)=0\qquad\textrm{as\:\ a\:\ function\:\ of}\; s_{0}\;\textrm{whenever}\;\mu\in\mathbb{Z_{\leq\textrm{0}}}\label{eq:zetaDerivSideB}\end{equation}

In {[}1{]}, section 4.2, we then performed explicit Cesaro calculations
of the root side, $r(s_{0},\mu)$, for $\zeta$ in the three cases
$\mu=0,-1$ and $-2$. In the first two cases we confirmed that 

\begin{equation}
r_{\zeta}(s_{0},0)=0\quad\textrm{and}\quad r_{\zeta}(s_{0},-1)=0\label{eq:zetaRootSideA}\end{equation}
also, at least modulo an estimate on $S(T)$ in the latter case which
holds conditional on RH; so that $\zeta$ continues to satisfy the
generalised root identities in these two cases. When $\mu=-2$, however,
we found that 

\begin{equation}
r_{\zeta}(s_{0},-2)=-\frac{1}{2}\neq0\label{eq:zetaRootSideB}\end{equation}
based on the assumption of the Riemann hypothesis; from which we concluded
that the Riemann hypothesis appears to be false. At no stage in {[}1{]},
however, did we consider the root side, $r_{\zeta}(s_{0},\mu)$, for
$\zeta$ for arbitrary $\mu\notin\mathbb{Z}$, which we noted represented
a gap in rigour in our arguments, since care must be taken in ascribing
Cesaro limits at isolated values of a parameter ($\mu$) in which
analytic continuation is being undertaken.

We thus now consider the possibility of arbitrary $\mu\notin\mathbb{Z}$
in the root identities for $\zeta$. In (\ref{eq:zetaDerivSideA})
this extension requires no further work for the derivative side, which
is convergent for arbitrary $\mu$ when $\Re(s_{0})>1$.

On the root side, we write $r_{\zeta}(s_{0},\mu)$ as the sum of three
separate contributions from $T$ (trivial roots), $P$ (simple pole)
and $NT$ (non-trivial roots) as in {[}1{]}:

\begin{equation}
r_{\zeta}(s_{0},\mu)=r_{T}(s_{0},\mu)+r_{P}(s_{0},\mu)+r_{NT}(s_{0},\mu)\label{eq:zetaRootSideC1}\end{equation}
where

\begin{eqnarray}
r_{T}(s_{0},\mu) & = & \textrm{e}^{i\pi\mu}\sum_{\{s_{0}-T\}}\frac{1}{(s_{0}-r_{i})^{\mu}}\:,\label{eq:zetaRootSideC2}\\
r_{P}(s_{0},\mu) & = & -\textrm{e}^{i\pi\mu}\frac{1}{(s_{0}-1)^{\mu}}\:,\:\textrm{and}\label{eq:zetaRootSideC3}\\
r_{NT}(s_{0},\mu) & = & \textrm{e}^{i\pi\mu}\sum_{\{s_{0}-NT\}}\frac{1}{(s_{0}-\rho_{i})^{\mu}}\label{eq:zetaRootSideC4}\end{eqnarray}

When $\Re(\mu)>1$, this is also classically convergent for generic
$s_{0}\in\mathbb{C}$ and so may be evaluated to any required degree
of accuracy by including sufficiently many roots, both trivial and
non-trivial. Here, by generic $s_{0}$ we mean $s_{0}$ not on either
$(-\infty,1)$ nor on any of the rays running from a non-trivial root
$\rho_{i}=\beta_{i}+i\gamma_{i}$ horizontally to the left towards
$-\infty$; these being branch cuts of terms in the expressions above.
Since it is all that is required for our arguments, and simultaneously
avoids all of these branch-cut issues, from now on we shall always
restrict consideration to $s_{0}$ in the half-plane $\Re(s_{0})>1$.

\subsection{The Case of $\Re(\mu)>1$}

For such $s_{0}$ we may then readily test numerically whether $\zeta$
does satisfy the generalised root identities also for arbitrary $\Re(\mu)>1$
by simply evaluating the expressions (\ref{eq:zetaDerivSideA}) and
(\ref{eq:zetaRootSideC1}) to the required degree of accuracy and
checking whether we do indeed have

\begin{equation}
d_{\zeta}(s_{0},\mu)=r_{\zeta}(s_{0},\mu)\label{eq:zetaGenRootIdentity}\end{equation}
For the derivative side this requires including sufficiently many
primes $p$ (the larger $s_{0}$ the fewer primes required to reach
a given accuracy); for the root side it requires including sufficiently
many trivial roots for $r_{T}$ and sufficiently many non-trivial
roots for $r_{NT}$ (the more positive $\mu>1$ is, the fewer such
roots required for a given level of accuracy).

In order to implement such numerical tests a list of primes and a
list of non-trivial roots must be used. Such lists can readily be
obtained and downloaded - e.g. a list of the first $2,001,052\: NT$
zeros can be downloaded as a text file from {[}2{]}, and in fact {[}3{]}
contains files listing up to the first $1,000,000,000\: NT$ zeros.

Using these we may then readily code up checks for random values of
$s_{0}>1$ and $\mu>1$ and see that $\zeta$ does indeed seem to
satisfy the generalised root identities (\ref{eq:zetaGenRootIdentity})
for arbitrary $s_{0},\mu$ in this region (starting with real $s_{0},\mu$
and then extending by analytical continuation).

Indeed Dana Pascovici's R-script for $\zeta$, available with this
paper, contains such code (currently using the first $10,000$ primes
(for $d_{\zeta}$) and the first $2,000,000$ trivial and $NT$ zeros
(for $r_{\zeta}$) but readily adaptable to use more of either if
desired). This thereby allows rapid, systematic checking and on the
basis of this we can confirm that $\zeta$ does indeed appear to satisfy
the generalised root identities for arbitrary $\Re(s_{0})>1$ when
$\Re(\mu)>1$, albeit that we cannot extend the testing down to $\mu$-values
too close to $1$ without radically increasing the number of roots
used in our numerical approximations for $r_{\zeta}(s_{0},\mu)$,
for obvious convergence reasons. The following pictures, for example,
illustrate the agreement of both real and imaginary parts of $d_{\zeta}(s_{0},\mu)$
and $r_{\zeta}(s_{0},\mu)$ for $s_{0}=5.1238$ for a range of $\mu$-values
greater than $1.5$.

\includegraphics[scale=0.5]{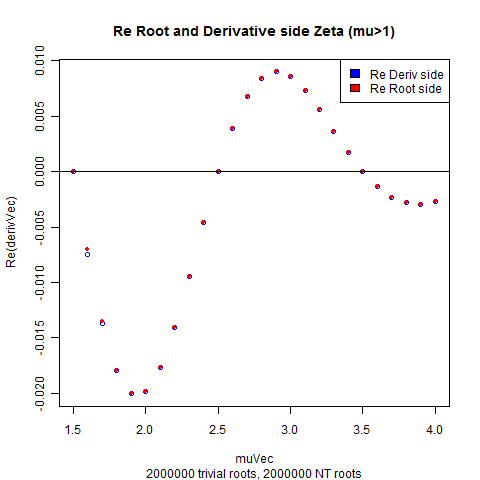}

\includegraphics[scale=0.5]{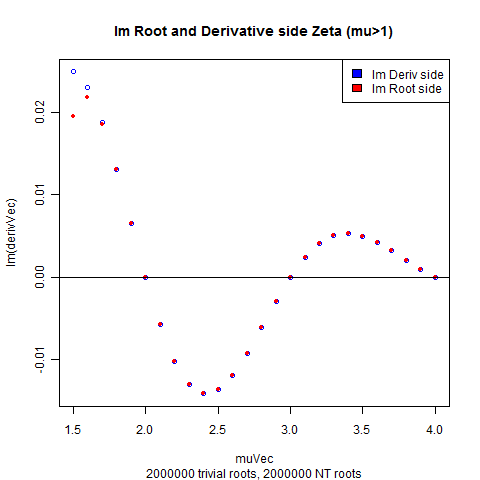}

\subsection{The Case of $\Re(\mu)<1$}

When $\Re(\mu)<1$ the evaluation of both $r_{T}$ and $r_{NT}$ on
the root side requires Cesaro methods, since the sums in (\ref{eq:zetaRootSideC2})
and (\ref{eq:zetaRootSideC4}) are no longer classically convergent.

\subsubsection{Formulae for numerical evaluation of $r_{\zeta}(s_{0},\mu)$ when
$\mu\in(-1,1)\setminus\{0\}$}

We take the components of the root side in turn:

\paragraph{Part (a) - $r_{T}(s_{0},\mu)$:}

For $r_{T}(s_{0},\mu)$ we use the Euler-McLaurin sum formula for
the partial sums $s_{T}$. Since the shifted roots occur at $s_{0}+2j$,
$j\in\mathbb{Z_{\geq\textrm{1}}}$, we let $z=s_{0}+2k+\alpha$ with
$0\leq\alpha<2$ and write

\begin{equation}
r_{T}(s_{0},\mu)=\textrm{e}^{i\pi\mu}\cdot\underset{z\rightarrow\infty}{Clim}\, s_{T}(s_{0},\mu;z)\label{eq:zetaRootSideD1}\end{equation}
where

\begin{equation}
s_{T}(s_{0},\mu;z)=\sum_{j=1}^{k}\frac{1}{(s_{0}+2j)^{\mu}}\label{eq:zetaRootSideD2}\end{equation}
Now, by Euler-McLaurin,

\begin{equation}
\sum_{j=1}^{k}\frac{1}{(s_{0}+2j)^{\mu}}=\left\{ \begin{array}{cc}
\frac{1}{2}\frac{(s_{0}+2k)^{1-\mu}}{1-\mu}+C(s_{0},\mu)+\frac{1}{2}(s_{0}+2k)^{-\mu}\\
\\-\frac{1}{6}\mu(s_{0}+2k)^{-\mu-1}\\
\\+\frac{\mu(\mu+1)(\mu+2)}{90}(s_{0}+2k)^{-\mu-3}-\ldots\end{array}\right\} \label{eq:zetaRootSideD3}\end{equation}
and so we can obtain the Cesaro limit in (\ref{eq:zetaRootSideD1})
by the usual Cesaro-style calculations on terms of the form $(s_{0}+2k+\alpha)^{-\mu-l}$
(see e.g. {[}1{]} or the earlier working for $\Gamma(z+1)$).

For example, considering just the case $0<\Re(\mu)<1$ initially,
here we have

\begin{eqnarray}
(s_{0}+2k+\alpha)^{1-\mu} & = & (s_{0}+2k)^{1-\mu}\left\{ 1+(1-\mu)\frac{\alpha}{(s_{0}+2k)}+\ldots\right\} \nonumber \\
\nonumber \\ & = & (s_{0}+2k)^{1-\mu}+o(1)\label{eq:zetaRootSideD4}\end{eqnarray}
and thus since pure powers $z^{\rho}\;(\rho\neq0)$ have generalised
Cesaro limit $0$, so

\begin{equation}
r_{T}(s_{0},\mu)=\textrm{e}^{i\pi\mu}\underset{k\rightarrow\infty}{lim}\left\{ \sum_{j=1}^{k}\frac{1}{(s_{0}+2j)^{\mu}}-\frac{1}{2}\frac{(s_{0}+2k)^{1-\mu}}{1-\mu}\right\} \quad(0<\Re(\mu)<1)\label{eq:zetaRootSideD5}\end{equation}
where this is now a classical limit.

Similar but more involved calculations can be performed on the lower
order terms to extend (\ref{eq:zetaRootSideD5}) to a formula applicable
for $-1<\Re(\mu)<0$, then $-2<\Re(\mu)<-1$ and so on. We shall only
consider the next case of $-1<\Re(\mu)<0$ here and in this case the
formula becomes

\begin{equation}
r_{T}(s_{0},\mu)=\textrm{e}^{i\pi\mu}\underset{k\rightarrow\infty}{lim}\left\{ \sum_{j=1}^{k}\frac{1}{(s_{0}+2j)^{\mu}}-\frac{1}{2}\frac{(s_{0}+2k)^{1-\mu}}{1-\mu}-\frac{1}{2}(s_{0}+2k)^{-\mu}\right\} \label{eq:zetaRootSideD6}\end{equation}

\paragraph{Part (b) - $r_{P}(s_{0},\mu)$:}

This is clearly given by 

\begin{equation}
r_{P}(s_{0},\mu)=-\textrm{e}^{i\pi\mu}\frac{1}{(s_{0}-1)^{\mu}}\label{eq:zetaRootSideE1}\end{equation}
for any $s_{0}$ and $\mu$.

\paragraph{Part (c) - $r_{NT}(s_{0},\mu)$:}

Here we use the Riemann-von Mangoldt formula for the counting function
$N(T)$ and mimic the argumentation in {[}1{]}. Write $NT=NT_{+}\cup NT_{-}$where
$NT_{+}$ refers to the subset of non-trivial zeros with positive
imaginary part and $NT_{-}$ to the subset with negative imaginary
part. We work initially on $NT_{+}$ and then simply state the corresponding
formulae for $NT_{-}$ which are derived in analogous fashion. For
a root $\rho_{i}=\beta_{i}+i\gamma_{i}\in NT_{+}$ we have, on writing
$\beta_{i}=\frac{1}{2}+\epsilon_{i}$, that

\begin{eqnarray}
\frac{1}{(s_{0}-(\beta_{i}+\gamma_{i}))^{\mu}} & = & (s_{0}-\frac{1}{2}-\epsilon_{i}-i\gamma_{i})^{-\mu}\nonumber \\
\nonumber \\ & = & \textrm{e}^{i\frac{\pi}{2}\mu}\gamma_{i}^{-\mu}\left(1+i\frac{(s_{0}-\frac{1}{2}-\epsilon_{i})}{\gamma_{i}}\right)^{-\mu}\nonumber \\
\nonumber \\ & = & \textrm{e}^{i\frac{\pi}{2}\mu}\gamma_{i}^{-\mu}\left\{ \begin{array}{cc}
1-i\mu\frac{(s_{0}-\frac{1}{2}-\epsilon_{i})}{\gamma_{i}}\\
\\-\frac{\mu(\mu+1)}{2!}\frac{(s_{0}-\frac{1}{2}-\epsilon_{i})^{2}}{\gamma_{i}^{2}}\\
\\+i\frac{\mu(\mu+1)(\mu+2)}{3!}\frac{(s_{0}-\frac{1}{2}-\epsilon_{i})^{3}}{\gamma_{i}^{3}}+\ldots\end{array}\right\} \label{eq:zetaRootSideF1}\end{eqnarray}
and therefore, bearing in mind that any roots off the critical line
occur in symmetric pairs either side of it and arguing as in {[}1{]},
we have

\begin{eqnarray}
\sum_{\{s_{0}-NT_{+}\}}\frac{M_{i}}{(s_{0}-\rho_{i})^{\mu}} & = & \textrm{e}^{i\frac{\pi}{2}\mu}\left\{ \begin{array}{cc}
\sum_{\{s_{0}-NT_{+}\}}M_{i}\gamma_{i}^{-\mu}\\
\\-i\mu(s_{0}-\frac{1}{2})\sum_{\{s_{0}-NT_{+}\}}M_{i}\gamma_{i}^{-\mu-1}\\
\\-\frac{\mu(\mu+1)(s_{0}-\frac{1}{2})^{2}}{2!}\sum_{\{s_{0}-NT_{+}\}}M_{i}\gamma_{i}^{-\mu-2}\\
\\-\frac{\mu(\mu+1)}{2!}\sum_{\{s_{0}-NT_{+}\}}M_{i}\gamma_{i}^{-\mu-2}\epsilon_{i}^{2}\\
\\+i\frac{\mu(\mu+1)(\mu+2)(s_{0}-\frac{1}{2})^{3}}{3!}\sum_{\{s_{0}-NT_{+}\}}M_{i}\gamma_{i}^{-\mu-3}\\
\\+i\frac{\mu(\mu+1)(\mu+2)(s_{0}-\frac{1}{2})}{2}\sum_{\{s_{0}-NT_{+}\}}M_{i}\gamma_{i}^{-\mu-3}\epsilon_{i}^{2}+\ldots\end{array}\right\} \nonumber \\
\label{eq:zetaRootSideF2}\end{eqnarray}

We next need to express the partial sums for these series in terms
of the parameter $T$ using the Riemann-von Mangoldt formula that 

\begin{equation}
N(T)=\check{N}(T)+S(T)+\frac{1}{\pi}\delta(T)\label{eq:zetaRootSideF3}\end{equation}
where, as in {[}1{]},

\begin{equation}
\check{N}(T)=\frac{T}{2\pi}\ln(\frac{T}{2\pi})-\frac{T}{2\pi}+\frac{7}{8}\label{eq:zetaRootSideF4}\end{equation}
and $S(T)$ is the argument of the zeta function, and

\begin{eqnarray}
\delta(T) & = & \frac{T}{4}\ln\left(1+\frac{1}{4T^{2}}\right)+\frac{1}{4}\tan^{-1}\left(\frac{1}{2T}\right)-\frac{T}{2}\int_{0}^{\infty}\frac{(\frac{1}{2}-\{u\})}{(u+\frac{1}{4})^{2}+(\frac{T}{2})^{2}}\textrm{d}u\nonumber \\
\label{eq:zetaRootSideF5}\\ & = & \frac{1}{48}\frac{1}{T}+\frac{7}{5760}\frac{1}{T^{3}}+\frac{31}{80640}\frac{1}{T^{5}}+o\left(\frac{1}{T^{5}}\right)\label{eq:zetaRootSideF6}\end{eqnarray}
Once we have such expressions we will then re-express them in terms
of the geometric variable z (as in {[}1{]}) and try to evaluate their
generalised Cesaro limits by removing eigenfunctions $z^{\rho},\: z^{\rho}\ln z$
etc and averaging the remainder (in $T$).

Note, however, that in combining with $N(T)$ only the sums in (\ref{eq:zetaRootSideF2})
that have no explicit $\epsilon_{i}$-dependence will be able to be
handled. This is because we certainly have no understanding of the
functional dependence of any $\epsilon_{i}(T)$ on $T$. In part for
this reason, and also for reasons of numerical tractability, for this
paper we shall only consider the extension of numerical checks to
$-1<\Re(\mu)<1$, $\mu\notin\mathbb{Z}$.

So now suppose $-1<\Re(\mu)<1$. Then in (\ref{eq:zetaRootSideF2})
we have only the sums $\sum_{\{s_{0}-NT_{+}\}}M_{i}\gamma_{i}^{-\mu}$
and $\sum_{\{s_{0}-NT_{+}\}}M_{i}\gamma_{i}^{-\mu-1}$ being divergent,
the former for all such $\mu$, the latter only when $\Re(\mu)<0$.
Take the case of $\sum_{\{s_{0}-NT_{+}\}}M_{i}\gamma_{i}^{-\mu}$
first. Considering just $NT_{+}$ initially, this generates a contribution
up to imaginary part $T$ of

\begin{eqnarray}
\int_{0}^{T}t^{-\mu}dN(t) & = & T^{-\mu}N(T)+\mu\int_{0}^{T}t^{-\mu-1}N(t)\:\textrm{d}t\label{eq:zetaRootSideG1}\\
\nonumber \\ & = & \left\{ \begin{array}{cc}
T^{-\mu}\check{N}(T)+\mu\int_{0}^{T}t^{-\mu-1}\check{N}(t)\:\textrm{d}t\\
\\+T^{-\mu}S(T)+\mu\int_{0}^{T}t^{-\mu-1}S(t)\:\textrm{d}t\\
\\+\frac{1}{\pi}T^{-\mu}\delta(T)+\frac{1}{\pi}\mu\int_{0}^{T}t^{-\mu-1}\delta(t)\:\textrm{d}t\end{array}\right\} \label{eq:zetaRootSideG2}\end{eqnarray}
Since $\delta(T)=O(\frac{1}{T})$, the terms involving $\delta(T)$
are classically convergent when $\Re(\mu)>-1$. And by arguments analogous
to those in {[}1{]} the terms involving $S(T)$ are also convergent
under some suitable power of the Cesaro averaging operator $P$ (they
are classically convergent when $\Re(\mu)>0$ and require one application
of $P$ when $-1<\Re(\mu)<0$).

The only terms where there are explicit eigenfunction divergences
that require removing are therefore those involving $\check{N}(T)$.
Here, from (\ref{eq:zetaRootSideF4}) we have, on writing $u=\frac{T}{2\pi}$,
using integration by parts \& changing variables to $v=\frac{t}{2\pi}$,
that

\[
T^{-\mu}\check{N}(T)+\mu\int_{0}^{T}t^{-\mu-1}\check{N}(t)\,\textrm{d}t=(2\pi)^{-\mu}\left\{ \begin{array}{cc}
\left[u^{1-\mu}\ln u-u^{1-\mu}+\frac{7}{8}u^{-\mu}\right]\\
+\mu\intop_{0}^{u}\left[\begin{array}{cc}
v^{-\mu}\ln v-v^{-\mu}\\
+\frac{7}{8}v^{-\mu-1}\end{array}\right]\textrm{d}v\end{array}\right\} \]

\begin{eqnarray}
 &  & =(2\pi)^{-\mu}\left\{ \begin{array}{cc}
\frac{1}{1-\mu}\left(\frac{T}{2\pi}\right)^{1-\mu}\ln\left(\frac{T}{2\pi}\right)\\
-\frac{1}{(1-\mu)^{2}}\left(\frac{T}{2\pi}\right)^{1-\mu}\end{array}\right\} \nonumber \\
\nonumber \\ &  & +C+o(1)\label{eq:zetaRootSideG3}\end{eqnarray}

Similarly, for $\sum_{\{s_{0}-NT_{+}\}}\gamma_{i}^{-\mu-1}$ we have
contribution

\begin{equation}
\int_{0}^{T}t^{-\mu-1}dN(t)=T^{-\mu-1}N(T)+(\mu+1)\int_{0}^{T}t^{-\mu-2}N(t)\:\textrm{d}t\label{eq:zetaRootSideG4}\end{equation}
and the only divergent terms come from the $\check{N}(T)$ terms,
namely

\begin{eqnarray}
T^{-\mu-1}\check{N}(T)+(\mu+1)\int_{0}^{T}t^{-\mu-2}\check{N}(t)\,\textrm{d}t & =(2\pi)^{-\mu-1} & \left\{ \begin{array}{cc}
-\frac{1}{\mu}\left(\frac{T}{2\pi}\right)^{-\mu}\ln\left(\frac{T}{2\pi}\right)\\
-\frac{1}{\mu^{2}}\left(\frac{T}{2\pi}\right)^{-\mu}\end{array}\right\} \nonumber \\
\nonumber \\ &  & +C+o(1)\label{eq:zetaRootSideG5}\end{eqnarray}

It follows from combining these calculations in (\ref{eq:zetaRootSideF2})
that, when $-1<\Re(\mu)<1$, the divergent part of $\sum_{\{s_{0}-NT_{+}\}}\frac{M_{i}}{(s_{0}-\rho_{i})^{\mu}}$
is given by $divgt_{NT_{+}}(s_{0},\mu;T)$ where

\begin{eqnarray}
divgt_{NT_{+}}(s_{0},\mu;T) & = & \textrm{e}^{i\frac{\pi}{2}\mu}(2\pi)^{-\mu}\left\{ \begin{array}{cc}
\frac{1}{1-\mu}\left(\frac{T}{2\pi}\right)^{1-\mu}\ln\left(\frac{T}{2\pi}\right)\\
-\frac{1}{(1-\mu)^{2}}\left(\frac{T}{2\pi}\right)^{1-\mu}\end{array}\right\} \nonumber \\
\nonumber \\ &  & +i\textrm{e}^{i\frac{\pi}{2}\mu}(2\pi)^{-\mu-1}(s_{0}-\frac{1}{2})\left\{ \begin{array}{cc}
\left(\frac{T}{2\pi}\right)^{-\mu}\ln\left(\frac{T}{2\pi}\right)\\
+\frac{1}{\mu}\left(\frac{T}{2\pi}\right)^{-\mu}\end{array}\right\} \nonumber \\
\label{eq:zetaRootSideG6}\end{eqnarray}

Since $\mu\notin\mathbb{Z}$ it is then clear that $divgt_{NT_{+}}(s_{0},\mu;T)$
can be expressed exactly in terms of powers of the geometric variable
$z=(s_{0}-\frac{1}{2})-iT$. In fact

\begin{eqnarray}
divgt_{NT_{+}}(s_{0},\mu;z) & = & \left\{ \begin{array}{cc}
i\,\frac{(2\pi)^{-\mu}}{1-\mu}\left(\frac{z}{2\pi}\right)^{1-\mu}\ln\left(\frac{z}{2\pi}\right)\\
\\-(2\pi)^{-\mu}\left\{ \frac{\frac{\pi}{2}}{1-\mu}+\frac{i}{(1-\mu)^{2}}\right\} \left(\frac{z}{2\pi}\right)^{1-\mu}\\
\\+(2\pi)^{-\mu-1}\frac{i}{\mu}(s_{0}-\frac{1}{2})\left(\frac{z}{2\pi}\right)^{-\mu}\end{array}\right\} \label{eq:zetaRootSideG7}\end{eqnarray}
and if $\mu\neq0$ in $(-1,1)$ all these terms have generalised Cesaro
limit $0$ so that 

\begin{equation}
\underset{z\rightarrow\infty}{Clim}\, divgt_{NT_{+}}(s_{0},\mu;z)=0\quad\quad\textrm{when}\;\mu\in(-1,1)\setminus\{0\}\label{eq:zetaRootSideG8}\end{equation}
In other words, for $\mu\in(-1,1)\setminus\{0\}$

\begin{equation}
\sum_{\{s_{0}-NT_{+}\}}\frac{M_{i}}{(s_{0}-\rho_{i})^{\mu}}=\underset{T\rightarrow\infty}{Clim}\left\{ \sum_{\Im(\rho_{i})<T}\,\frac{M_{i}}{(s_{0}-\rho_{i})^{\mu}}-divgt_{NT_{+}}(s_{0},\mu;T)\right\} \label{eq:zetaRootSideG9}\end{equation}
where $divgt_{NT_{+}}(s_{0},\mu;T)$ is as per (\ref{eq:zetaRootSideG6})
and this is in fact a classical limit if $\mu>0$ and requires only
one application of the averaging operator, $P_{T}$, if $-1<\mu<0$.

In similar fashion, for the $NT$ roots below the real axis, $NT_{_{-}}$,
we have

\begin{eqnarray}
\sum_{\{s_{0}-NT_{-}\}}\frac{M_{i}}{(s_{0}-\rho_{i})^{\mu}} & = & \textrm{e}^{-i\frac{\pi}{2}\mu}\left\{ \begin{array}{cc}
\sum_{\{s_{0}-NT_{-}\}}\gamma_{i}^{-\mu}\\
\\+i\mu(s_{0}-\frac{1}{2})\sum_{\{s_{0}-NT_{-}\}}\gamma_{i}^{-\mu-1}\\
\\-\frac{\mu(\mu+1)}{2!}(s_{0}-\frac{1}{2})^{2}\sum_{\{s_{0}-NT_{-}\}}\gamma_{i}^{-\mu-2}\\
\\-\frac{\mu(\mu+1)}{2!}\sum_{\{s_{0}-NT_{-}\}}\gamma_{i}^{-\mu-2}\epsilon_{i}^{2}\\
\\-i\frac{\mu(\mu+1)(\mu+2)}{3!}(s_{0}-\frac{1}{2})^{3}\sum_{\{s_{0}-NT_{-}\}}\gamma_{i}^{-\mu-3}\\
\\-i\frac{\mu(\mu+1)(\mu+2)}{2}(s_{0}-\frac{1}{2})\sum_{\{s_{0}-NT_{-}\}}\gamma_{i}^{-\mu-3}\epsilon_{i}^{2}+\ldots\end{array}\right\} \nonumber \\
\label{eq:zetaRootSideG10}\end{eqnarray}
and thus, by identical reasoning, for $\mu\in(-1,1)\setminus\{0\}$
the divergent piece of the partial sums to {}``height'' $\tilde{T}$
for $NT_{-}$ is given by

\begin{eqnarray}
divgt_{NT_{-}}(s_{0},\mu;\tilde{T}) & = & \textrm{e}^{-i\frac{\pi}{2}\mu}(2\pi)^{-\mu}\left\{ \begin{array}{cc}
\frac{1}{1-\mu}\left(\frac{\tilde{T}}{2\pi}\right)^{1-\mu}\ln\left(\frac{\tilde{T}}{2\pi}\right)\\
-\frac{1}{(1-\mu)^{2}}\left(\frac{\tilde{T}}{2\pi}\right)^{1-\mu}\end{array}\right\} \nonumber \\
\nonumber \\ &  & -i\textrm{e}^{-i\frac{\pi}{2}\mu}(2\pi)^{-\mu-1}(s_{0}-\frac{1}{2})\left\{ \begin{array}{cc}
\left(\frac{\tilde{T}}{2\pi}\right)^{-\mu}\ln\left(\frac{\tilde{T}}{2\pi}\right)\\
+\frac{1}{\mu}\left(\frac{\tilde{T}}{2\pi}\right)^{-\mu}\end{array}\right\} \nonumber \\
\label{eq:zetaRootSideG11}\end{eqnarray}
For $\mu\neq0$ this can be expressed exactly in terms of the geometric
variable $\tilde{z}=(s_{0}-\frac{1}{2})+i\tilde{T}$ as

\begin{equation}
divgt_{NT_{-}}(s_{0},\mu;\tilde{z})=\left\{ \begin{array}{cc}
-i\,\frac{(2\pi)^{-\mu}}{1-\mu}\left(\frac{\tilde{z}}{2\pi}\right)^{1-\mu}\ln\left(\frac{\tilde{z}}{2\pi}\right)\\
\\-(2\pi)^{-\mu}\left\{ \frac{\frac{\pi}{2}}{1-\mu}-\frac{i}{(1-\mu)^{2}}\right\} \left(\frac{\tilde{z}}{2\pi}\right)^{1-\mu}\\
\\-(2\pi)^{-\mu-1}\frac{i}{\mu}(s_{0}-\frac{1}{2})\left(\frac{\tilde{z}}{2\pi}\right)^{-\mu}\end{array}\right\} \label{eq:zetaRootSideG12}\end{equation}
with generalised Cesaro limit $0$, i.e.

\begin{equation}
\underset{\tilde{z}\rightarrow\infty}{Clim}\, divgt_{NT_{-}}(s_{0},\mu;\tilde{z})=0\quad\quad\textrm{when}\;\mu\in(-1,1)\setminus\{0\}\label{eq:zetaRootSideG13}\end{equation}
so that, for $\mu\in(-1,1)\setminus\{0\}$

\begin{equation}
\sum_{\{s_{0}-NT_{-}\}}\frac{M_{i}}{(s_{0}-\rho_{i})^{\mu}}=\underset{\tilde{T}\rightarrow\infty}{Clim}\left\{ \sum_{\Im(\rho_{i})>-\tilde{T}}\,\frac{M_{i}}{(s_{0}-\rho_{i})^{\mu}}-divgt_{NT_{-}}(s_{0},\mu;\tilde{T})\right\} \label{eq:zetaRootSideG14}\end{equation}
where $divgt_{NT_{-}}(s_{0},\mu;\tilde{T})$ is as given in (\ref{eq:zetaRootSideG11})
and again this is a classical limit if $\mu>0$ and requires one application
of the averaging operator, $P_{\tilde{T}}$, if $-1<\mu<0$.

If we write $\tilde{r}_{NT_{+}}(s_{0},\mu)$ and $\tilde{r}_{NT_{-}}(s_{0},\mu)$
for the expressions on the right in (\ref{eq:zetaRootSideG9}) and
(\ref{eq:zetaRootSideG14}) we thus have

\begin{equation}
r_{NT}(s_{0},\mu)=\textrm{e}^{i\pi\mu}\sum_{\{s_{0}-NT\}}\frac{M_{i}}{(s_{0}-\rho_{i})^{\mu}}=r_{NT_{+}}(s_{0},\mu)+r_{NT_{-}}(s_{0},\mu)\label{eq:zetaRootSideG15}\end{equation}
where

\begin{equation}
r_{NT_{+}}(s_{0},\mu)=\textrm{e}^{i\pi\mu}\tilde{r}_{NT_{+}}(s_{0},\mu)\quad\quad\textrm{and}\quad\quad r_{NT_{-}}(s_{0},\mu)=\textrm{e}^{i\pi\mu}\tilde{r}_{NT_{-}}(s_{0},\mu)\label{eq:zetaRootSideG16}\end{equation}
and the expressions in (\ref{eq:zetaRootSideG9}) and (\ref{eq:zetaRootSideG14})
then allow numerical computation of $r_{NT}(s_{0},\mu)$ also for
$\mu\in(-1,1)\setminus\{0\}$.

\subsubsection{Results of numerical tests for $\zeta$ for $\mu\in(-1,1)\setminus\{0\}$}

Based on the results in parts (a)-(c) in section 4.2.1 we have proceeded
to implement code for numerical evaluation of $r_{\zeta}(s_{0},\mu)$
and thus for testing of the generalised root identities for $\zeta$
for $\mu\in(-1,1)\setminus\{0\}$. The following general observations
apply to these efforts:

(a) With only $2,000,000$ $NT_{+}$ roots, convergence is too slow
to obtain reasonable approximate values on the root side when $\mu$
approaches $0$ from \textit{above}. This is for the same reason that
we couldn't obtain reasonable approximations even in the convergent
region ($\mu>1$) when $\mu$ approached $1$ from above. When $\mu$
is just above an integer (e.g. $0$) the residual piece, even after
throwing away the divergent terms, converges too slowly (like $T^{-\mu}$
or $T^{-\mu}S(T)$ when $\mu>0$ small) to get reasonable convergence
with only $2,000,000$ $NT_{+}$ roots.

However, when $\mu$ approaches an integer (e.g. $0$) from \textit{below},
one can get accurate numerical results, since on throwing away the
extra divergent terms that appear when $\mu$ passes through the integer
value, the remaining residual term decays much more rapidly (like
$T^{-1-\mu}$ or $T^{-1-\mu}S(T)$ when $\mu<0$ small).

Thus for our numerical tests in this section we consider only $\mu$-regions
just to the \textit{left} of $\mu=1$ and $0$. This also helps get
accurate evaluations of the trivial-roots contributions, $r_{T}(s_{0},\mu)$,
with our choice of $2,000,000$ trivial roots.

(b) As noted, when $0<\mu<1$ we only need to remove divergent pieces
and then take classical limits (and in doing this we only subtract
those terms in the expressions (\ref{eq:zetaRootSideD6}), (\ref{eq:zetaRootSideG6})
and (\ref{eq:zetaRootSideG11}) which are themselves divergent when
$0<\mu<1$).

When $-1<\mu<0$, however, for the $NT$-root contributions we need
to remove the divergent pieces in (\ref{eq:zetaRootSideG6}) and (\ref{eq:zetaRootSideG11})
and then apply the Cesaro averaging operator (in $T$ or $\tilde{T}$)
once to the resulting partial-sum function on the critical line, whose
step-jumps on crossing roots are now bigger than $1$ and grow. Alternatively
this can be viewed as being needed in order to handle the oscillatory
contributions from the $S(T)$-components of these expressions.

For $\mu$ small negative and approaching $0$, however, we see that
the step-discontinuities approach height $1$ and the resulting partial-sum
functions become very close to piecewise linear functions with negative
slope between root-steps. This can be seen clearly, for example, in
the following picture:

\includegraphics[scale=0.32]{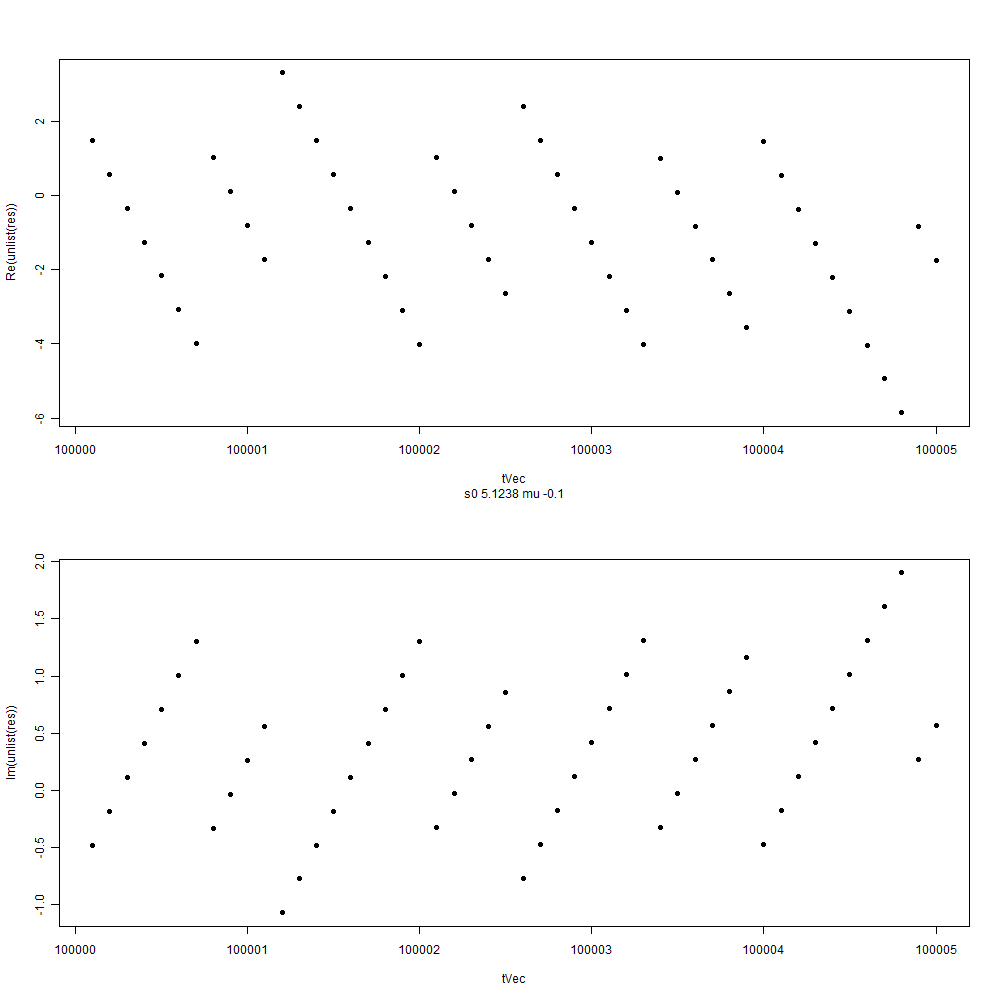}

This shows the partial-sum function for $\tilde{r}_{NT_{+}}(s_{0},\mu)$
when $s_{0}=5.1238$ and $\mu=-0.1$ on the interval $100,000<T<100,005$.
We see that by using the trapezoidal rule between each pair of successive
$NT$-roots we can readily implement code to perform the required
Cesaro averaging approximately for $\mu<0$ small, and hence numerically
evaluate the contributions from $r_{NT_{+}}(s_{0},\mu)$ and $r_{NT_{-}}(s_{0},\mu)$.
Using these we are then able to verify whether $\zeta$ does indeed
continue to satisfy the generalised root identities also when $\mu<0$
small and whether the root-side of these identities, $r_{\zeta}(s_{0},\mu)$,
does also tend to $0$ as $\mu\rightarrow0^{-}$, like the derivative
side $d_{\zeta}(s_{0},\mu)$.

(c) With the observations from (a) and (b) in mind, in implementing
our numerical root-side calculations to test the root-identities for
$\zeta$ we have:

(i) implemented this within the $\zeta$ R-script for $0<\mu<1$ and
in order to generate pictures of the partial-sum function in $T$
when $-1<\mu<0$ (as for example in the figure above), but

(ii) implemented the full calculations involving Cesaro averaging
when $-1<\mu<0$ in VBA code in XL2007; this can be adapted to the
more powerful framework of R, but is sufficient for our numerical
experiments in this paper and has the virtue of enhanced transparency
of the code.

Both the XL file ({}``RootIdentitiesZeta\_mu\_-1To0\_TestsB.xlsm'')
and the zeta R-script, along with their required source files, are
made available with this paper to allow a reader readily to examine
the code involved, replicate the claimed results and conduct further
independent tests if desired (see Appendix 5.1 for further brief instructions
on doing this).

\paragraph{Results:}

(a) When $0<\mu<1$ we get close agreement between $d_{\zeta}(s_{0},\mu)$
and $r_{\zeta}(s_{0},\mu)$ for $\mu>0.75$. For $\mu<0.75$, as noted,
$2,000,000$ $NT$ roots is insufficient to gain a good approximation
to the $NT$-root contributions to $r_{\zeta}(s_{0},\mu)$, with the
situation growing worse as $\mu\rightarrow0^{+}$. The following pictures
show how the convergence on the root side improves appreciably as
$\mu\rightarrow1^{-}$ and agreement between $d_{\zeta}(s_{0},\mu)$
and $r_{\zeta}(s_{0},\mu)$ thus becomes progressively better as $\mu\rightarrow1^{-}$:

\includegraphics[scale=0.5]{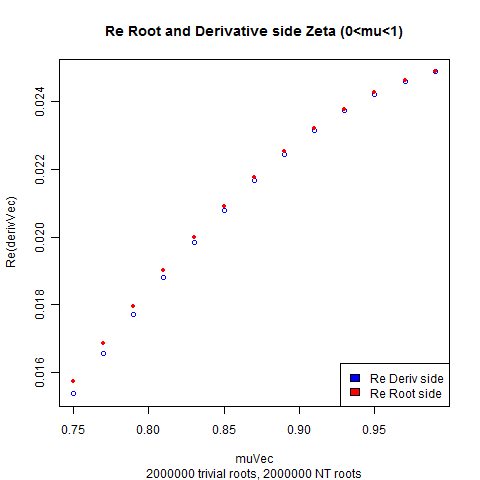}

\includegraphics[scale=0.5]{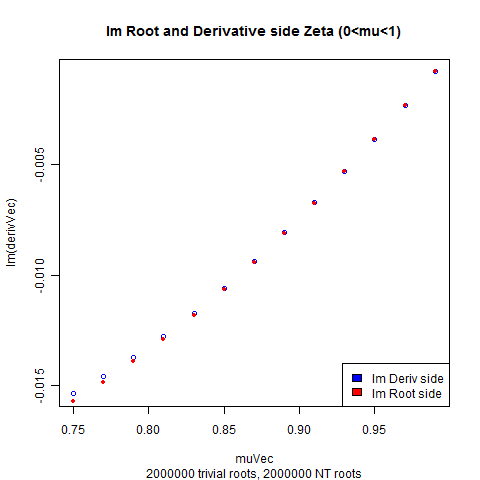}

We believe that these provide convincing evidence that $\zeta$ does
indeed satisfy the generalised root identities also for $0<\mu<1$;
this could be further tested, over a wider $\mu$-range, by using
more than $2,000,000$ trivial and $NT$-roots.

(b) When $-1<\mu<0$ we likewise get reasonably close agreement between
$d_{\zeta}(s_{0},\mu)$ and $r_{\zeta}(s_{0},\mu)$ for $\mu>-0.2$
after Cesaro averaging, but $2,000,000$ $NT$-roots is insufficient
to obtain reasonable root side values when $\mu<-0.2$: 

\includegraphics[scale=0.65]{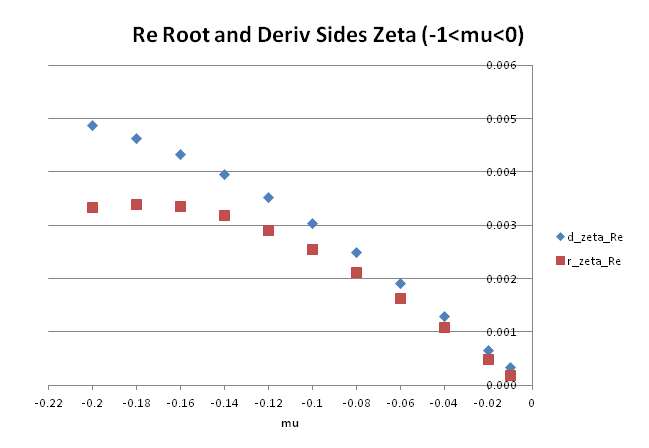}

\includegraphics[scale=0.65]{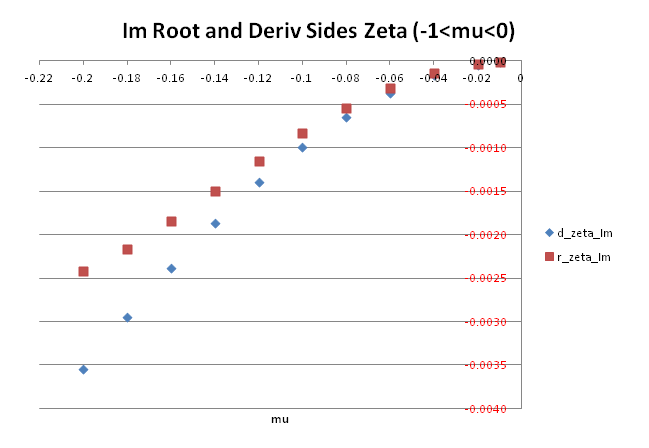}

These pictures are again for $s_{0}=5.1238$. Here, and in all testing
described for $-1<\mu<0$, the approximations are obtained using $10,000$
primes on the derivative side with a truncation threshold of $10$
(so that $d_{\zeta}(s_{0},\mu)$ values should be very accurate),
while on the root side we use $2,000,000$ trivial roots for $r_{T}$
and an average over the next $10,000$ roots starting at the $1,000,000^{th}$
$NT$-root, with 5 subintervals in the trapezoidal rule and an offset
of $10^{-9}$ (to avoid the step-jumps at roots) in the calculation
of $r_{NT_{+}}$ and $r_{NT_{-}}$.

The following tables show further detail on the breakdown of the calculations
for $d_{\zeta}(s_{0},\mu)$ and $r_{\zeta}(s_{0},\mu)$ when $\mu$
is very small, namely $\mu=-0.001$, $\mu=-0.0001$ and $\mu=-0.00001$:
(insert tables here)

\begin{flushleft}
\begin{tabular}{|c|c|c|c|c|c|c|}
\hline 
{\small $\mu$} & {\footnotesize $\Re(d_{\zeta})$} & {\footnotesize $\Re(r_{\zeta})$} & {\footnotesize $\Re(r_{T})$} & {\footnotesize $\Re(r_{P})$} & {\footnotesize $\Re(r_{NT_{+}})$} & {\footnotesize $\Re(r_{NT_{-}})$}\tabularnewline
\hline
\hline 
{\scriptsize -0.001} & {\scriptsize 0.00003300} & {\scriptsize -0.00009450} & {\scriptsize -3.064361} & {\scriptsize -1.001413} & {\scriptsize 4.348981 } & {\scriptsize -0.283302}\tabularnewline
\hline 
{\scriptsize -0.0001} & {\scriptsize 0.00000330} & {\scriptsize -0.00012270} & {\scriptsize -3.062147} & {\scriptsize -1.000142} & {\scriptsize 4.343429} & {\scriptsize -0.281263}\tabularnewline
\hline 
{\scriptsize -0.00001} & {\scriptsize 0.00000033 } & {\scriptsize -0.00012547 } & {\scriptsize -3.061925} & {\scriptsize -1.000014} & {\scriptsize 4.342870} & {\scriptsize -0.281056}\tabularnewline
\hline
\end{tabular}
\par\end{flushleft}

\begin{flushleft}
\begin{tabular}{|c|c|c|c|c|c|c|}
\hline 
{\small $\mu$} & {\footnotesize $\Im(d_{\zeta})$} & {\footnotesize $\Im(r_{\zeta})$} & {\footnotesize $\Im(r_{T})$} & {\footnotesize $\Im(r_{P})$} & {\footnotesize $\Im(r_{NT_{+}})$} & {\footnotesize $\Im(r_{NT_{-}})$}\tabularnewline
\hline
\hline 
{\scriptsize -0.001} & {\scriptsize -0.0000001037 } & {\scriptsize 0.0000002974} & {\scriptsize 0.009627 } & {\scriptsize 0.003146} & {\scriptsize 737.242} & {\scriptsize -737.255}\tabularnewline
\hline 
{\scriptsize -0.0001} & {\scriptsize -0.0000000010} & {\scriptsize 0.0000000385} & {\scriptsize 0.000962} & {\scriptsize 0.000314} & {\scriptsize 7,360.36} & {\scriptsize -7,360.36}\tabularnewline
\hline 
{\scriptsize -0.00001} & {\scriptsize -0.0000000000} & {\scriptsize 0.0000000039} & {\scriptsize 0.000096} & {\scriptsize 0.000031} & {\scriptsize 73,591.4} & {\scriptsize -73,591.4}\tabularnewline
\hline
\end{tabular}
\par\end{flushleft}

\begin{flushleft}
As $\mu\rightarrow0^{-}$ we see that:
\par\end{flushleft}

(i) Agreement between $d_{\zeta}(s_{0},\mu)$ and $r_{\zeta}(s_{0},\mu)$
becomes very close and both quantities approach $0$.%
\footnote{albeit there seems to be a threshold for the root side reflecting
the limits of accuracy we can obtain using only $2,000,000$ trivial
roots and an average over only $10,000$ roots starting at the $1,000,000^{th}$
$NT$-root in our approximate calculations%
}This provides strong evidence, first, that $\zeta$ does continue
to satisfy the generalised root identities for $\mu<0$ and, secondly,
that the values of $d_{\zeta}(s_{0},\mu)$ (and hence also $r_{\zeta}(s_{0},\mu)$)
behave continuously across the integer point $\mu=0$, with these
values approaching the limiting value of $0$ when $\mu=0$, as calculated
in section 4.1 in {[}1{]}.

(ii) Note in regard to (i) that while the root-side does approach
$0$ as $\mu\rightarrow0^{-}$, further confidence in the correctness
of the root identities comes from the fact that both the real and
imaginary parts of $r_{\zeta}(s_{0},\mu)$ are made up from the addition
of disparate elements ($r_{T}$, $r_{P}$, $r_{NT_{+}}$ and $r_{NT_{-}}$)
whose contributions (see the above tables) do not themselves become
small (and in some cases appear to diverge as $\mu\rightarrow0^{-}$).

(iii) In this regard note further that as $\mu$ goes from $-0.001$
to $-0.00001$ the value of $r_{T}$ more and more closely approaches
$-\frac{1}{2}s_{0}-\frac{1}{2}$ (as per equation (49) for $\mu=0$
in {[}1{]}); the value of $r_{P}$ more and more closely approaches
$-1$ (which it equals when $\mu=0$); and the value of $r_{NT}=r_{NT_{+}}+r_{NT_{-}}$
more and more closely approaches $\frac{1}{2}s_{0}+\frac{3}{2}$ (as
per equation (65) for $\mu=0$ in {[}1{]}). Thus the behaviour of
each component of the root side as $\mu\rightarrow0^{-}$ does appear
to continuously approach the behaviour calculated in {[}1{]} for these
components in the limiting case $\mu=0$.

\subsubsection{A conjecture regarding the appearance of the $\ln\left(\frac{z}{\tilde{z}}\right)$
divergences in the calculation of $r_{\zeta}(s_{0},\mu)$ when $\mu=0$
and numerical test of this conjecture}

In {[}1{]} the calculations of the $NT$-root contributions to the
root sides, $r_{\zeta}(s_{0},\mu)$, of the $\zeta$ root identities
at $\mu=0,-1\:\textrm{and}\:-2$ all relied on claiming that a value
of $0$ should be ascribed to the 2d Cesaro limit

\begin{equation}
\underset{z,\tilde{z}\rightarrow\infty}{Clim}\ln\left(\frac{z}{\tilde{z}}\right)=0\label{eq:2dCesaroLimit}\end{equation}
But in {[}1{]} this claim is only justified in a highly heuristic
fashion.%
\footnote{or claimed as an \textit{implication} of $\zeta$ satisfying the $\mu=0$
root identity, which can then at least be legitimately used for the
$\mu=-1$ and $\mu=-2$ calculations%
} By considering the case of non-integer $\mu$ as $\mu$ approaches
$0$, however, we may formulate a precise conjecture as to how the
expression $\ln\left(\frac{z}{\tilde{z}}\right)$ arises naturally
in the limit and why it would have to have generalised 2d Cesaro limit
$0$ in order to ensure continuity of the analytic continuation of
$r_{\zeta}(s_{0},\mu)$ in $\mu$ across $\mu=0$.

To see this, suppose $\mu\neq0$ is very small and consider the expressions
$\tilde{r}_{NT_{+}}(s_{0},\mu)$ and $\tilde{r}_{NT_{-}}(s_{0},\mu)$
defined by (\ref{eq:zetaRootSideG9}) and (\ref{eq:zetaRootSideG14}).
With $s_{0}$ real initially as usual, it is easy to see from the
expression for $divgt_{NT_{+}}(s_{0},\mu;T)$ in (\ref{eq:zetaRootSideG6})
and for $divgt_{NT_{-}}(s_{0},\mu;\tilde{T})$ in (\ref{eq:zetaRootSideG11})
that $\tilde{r}_{NT_{+}}(s_{0},\mu)$ and $\tilde{r}_{NT_{-}}(s_{0},\mu)$
must be complex conjugates:

\begin{equation}
\textrm{If}\quad\tilde{r}_{NT_{+}}(s_{0},\mu)=a+ib\qquad\textrm{then}\qquad\tilde{r}_{NT_{-}}(s_{0},\mu)=a-ib\label{eq:Conjec1}\end{equation}
Since the real part of $r_{NT}(s_{0},\mu)$ at $\mu=0$ is $\frac{1}{2}s_{0}+\frac{3}{2}$
(as shown in equation (63) in {[}1{]}) so, for the real parts, continuity
of analytic continuation across $\mu=0$ is equivalent to requiring
that

\[
\underset{\mu\rightarrow0}{lim}\,\Re(\textrm{e}^{i\pi\mu}(\tilde{r}_{NT_{+}}(s_{0},\mu)+\tilde{r}_{NT_{-}}(s_{0},\mu)))=\frac{1}{2}s_{0}+\frac{3}{2}\]
or, in light of (\ref{eq:Conjec1}) and the fact that $\underset{\mu\rightarrow0}{lim}\,\textrm{e}^{i\pi\mu}=1,$
that

\begin{equation}
\underset{\mu\rightarrow0}{lim}\,\Re(\tilde{r}_{NT_{+}}(s_{0},\mu))=\frac{1}{4}s_{0}+\frac{3}{4}\label{eq:Conjec2}\end{equation}

As for the imaginary part of $r_{NT}$, this automatically tends to
$0$ since by (\ref{eq:Conjec1}) the imaginary parts of $\tilde{r}_{NT_{+}}(s_{0},\mu)$
and $\tilde{r}_{NT_{-}}(s_{0},\mu)$ cancel and $\textrm{e}^{i\pi\mu}\rightarrow1$
as $\mu\rightarrow0$ as noted. Continuity of analytic continuation
across $\mu=0$ thus requires that the overall imaginary part of $r_{NT}$
also be $0$ at $\mu=0$, and since by equation (63) in {[}1{]} this
is precisely given by

\begin{equation}
\Im(r_{NT}(s_{0},0))=-\frac{1}{2\pi}(s_{0}-\frac{1}{2})\cdot\underset{z,\tilde{z}\rightarrow\infty}{Clim}\ln\left(\frac{z}{\tilde{z}}\right)\label{eq:Conjec3}\end{equation}
so this shows why we must have $\underset{z,\tilde{z}\rightarrow\infty}{Clim}\ln\left(\frac{z}{\tilde{z}}\right)=0$.

But we can go further and, by considering $\tilde{r}_{NT_{+}}(s_{0},\mu)$
and $\tilde{r}_{NT_{-}}(s_{0},\mu)$ separately we can conjecture
as to how the $\ln\left(\frac{z}{\tilde{z}}\right)$ term arises naturally
in the limit. Recall that (\ref{eq:zetaRootSideG9}) may be re-expressed
as saying that for $\mu\neq0$ small

\begin{equation}
\sum_{\Im(\rho_{i})<T}\,\frac{M_{i}}{(s_{0}-\rho_{i})^{\mu}}=divgt_{NT_{+}}(s_{0},\mu;z)+\tilde{r}_{NT_{+}}(s_{0},\mu)+R(s_{0},\mu;T)\label{eq:Conjec4}\end{equation}
where $divgt_{NT_{+}}(s_{0},\mu;z)$ is as given in (\ref{eq:zetaRootSideG7})
and where either $R=o(1)$ (if $0<\mu<1$) or $P[R](s_{0},\mu;T)\rightarrow0$
as $T\rightarrow\infty$ (if $-1<\mu<0$).

On the other hand, setting $u=\frac{z}{2\pi}$ and combining result
(60) from {[}1{]} and the further observation from the same section
of {[}1{]} that $P[S](T)\rightarrow0$ as $T\rightarrow\infty$, we
have the corresponding claim that at $\mu=0$

\begin{equation}
\sum_{\Im(\rho_{i})<T}\,\frac{M_{i}}{(s_{0}-\rho_{i})^{0}}=iu\cdot\ln u-(\frac{\pi}{2}+i)u-\frac{i}{2\pi}(s_{0}-\frac{1}{2})\ln u+(\frac{1}{4}s_{0}+\frac{3}{4})+R(s_{0},0;T)\label{eq:Conjec5}\end{equation}
where $P[R](s_{0},0;T)\rightarrow0$ as $T\rightarrow\infty$.

But now it is easy to see that, as $\mu\rightarrow0$, the first and
second terms in (\ref{eq:zetaRootSideG7}) turn into the corresponding
terms in (\ref{eq:Conjec5}) while we saw in (\ref{eq:Conjec2}) that
the real part of $\tilde{r}_{NT_{+}}(s_{0},\mu)$ must limit to the
term $\frac{1}{4}s_{0}+\frac{3}{4}$. As for the final term in (\ref{eq:zetaRootSideG7}),
using the the Taylor series expansion

\begin{equation}
\frac{1}{\mu}\, u^{-\mu}=\frac{1}{\mu}\left\{ 1-\mu\ln u+O(\mu^{2})\right\} =\frac{1}{\mu}-\ln u+O(\mu)\label{eq:Conjec6}\end{equation}
we see that, as $\mu\rightarrow0$, it formally turns into precisely
the term $-\frac{i}{2\pi}(s_{0}-\frac{1}{2})\ln u$ in (\ref{eq:Conjec5})
together with an additional pure imaginary constant term $\frac{i}{2\pi}(s_{0}-\frac{1}{2})\frac{1}{\mu}$
which blows up as $\mu\rightarrow0$.

In order to match (\ref{eq:Conjec5}) in the limit as $\mu\rightarrow0$
we would therefore have to have $\tilde{r}_{NT_{+}}(s_{0},\mu)$ containing
an exactly offsetting imaginary component, $-\frac{i}{2\pi}(s_{0}-\frac{1}{2})\frac{1}{\mu}$,
and so, combined with (\ref{eq:Conjec2}), we conjecture that

\paragraph{Conjecture:}

\textit{As} $\mu\rightarrow0$ 

\[
\tilde{r}_{NT_{+}}(s_{0},\mu)\sim(\frac{1}{4}s_{0}+\frac{3}{4})-\frac{i}{2\pi}(s_{0}-\frac{1}{2})\frac{1}{\mu}\]
\textit{meaning more precisely that}

\begin{equation}
\Re(\tilde{r}_{NT_{+}}(s_{0},\mu))=\frac{1}{4}s_{0}+\frac{3}{4}+o(1)\label{eq:Conjec7}\end{equation}
\textit{and}

\begin{equation}
\Im(\tilde{r}_{NT_{+}}(s_{0},\mu))=-\frac{i}{2\pi}(s_{0}-\frac{1}{2})\frac{1}{\mu}+C+o(1)\label{eq:Conjec8}\end{equation}
\textit{for some constant} $C$.

In this perspective then, the imaginary parts of $\tilde{r}_{NT_{+}}(s_{0},\mu)$
and $\tilde{r}_{NT_{-}}(s_{0},\mu)$ are both diverging as $\mu\rightarrow0$
but in exactly offsetting ways, and the residual $\ln\left(\frac{z}{\tilde{z}}\right)$
divergence which appears in handling the root side of the $\mu=0$
root identity is arising naturally from applying the Taylor series
expansion (\ref{eq:Conjec6}) to each of the final $u^{-\mu}$ and
$\tilde{u}^{-\mu}$ terms in the divergent pieces (\ref{eq:zetaRootSideG7})
and (\ref{eq:zetaRootSideG12}) for $\tilde{r}_{NT_{+}}(s_{0},\mu)$
and $\tilde{r}_{NT_{-}}(s_{0},\mu)$ and then taking the limit as
$\mu\rightarrow0$. Then $\ln\left(\frac{z}{\tilde{z}}\right)$ naturally
has to have 2d Cesaro limit $0$ since it arises from a limit in $\mu$
of the expression $u^{-\mu}-\tilde{u}^{-\mu}$ which clearly has Cesaro
limit $0$ at all $\mu$.%
\footnote{A formal L'hopital's calculation, as per section 4.3 of {[}4{]}, also
then recovers $\frac{1}{4}s_{0}+\frac{3}{4}$ in the limit as $\mu\rightarrow0$
under this conjecture%
}

\paragraph{Results:}

We tested this conjecture by calculating the values of $\Re(\tilde{r}_{NT_{+}}(s_{0},\mu))$
and $\Im(\tilde{r}_{NT_{+}}(s_{0},\mu))$ for $\mu=-0.001$, $\mu=-0.0001$
and $\mu=-0.00001$ as shown in the following tables (for $s_{0}=5.1238$
again): (insert tables here)

\begin{flushleft}
\begin{tabular}{|c|c|c|}
\hline 
{\scriptsize $\mu$} & {\scriptsize $\Re(\tilde{r}_{NT_{+}})$} & {\scriptsize Pred. Leading Order}\tabularnewline
\hline
\hline 
{\scriptsize -0.001} & {\scriptsize 2.0328495452 } & {\scriptsize 2.03095}\tabularnewline
\hline 
{\scriptsize -0.0001} & {\scriptsize 2.0310831931 } & {\scriptsize 2.03095}\tabularnewline
\hline 
{\scriptsize -0.00001} & {\scriptsize 2.0309067012 } & {\scriptsize 2.03095}\tabularnewline
\hline
\end{tabular}
\par\end{flushleft}

\begin{flushleft}
\begin{tabular}{|c|c|c|}
\hline 
{\scriptsize $\mu$} & {\scriptsize $\Im(\tilde{r}_{NT_{+}})$} & {\scriptsize Pred. Leading Order}\tabularnewline
\hline
\hline 
{\scriptsize -0.001} & {\scriptsize 745.20948403} & {\scriptsize 735.90063 }\tabularnewline
\hline 
{\scriptsize -0.0001} & {\scriptsize 7,361.1449162} & {\scriptsize 7,359.0063 }\tabularnewline
\hline 
{\scriptsize -0.00001} & {\scriptsize 73,591.493604} & {\scriptsize 73,590.063}\tabularnewline
\hline
\end{tabular}
\par\end{flushleft}

We clearly see that both parts of the above conjecture appear to be
confirmed numerically with $\Re(\tilde{r}_{NT_{+}}(s_{0},\mu))$ becoming
progressively closer and closer to $\frac{1}{4}s_{0}+\frac{3}{4}=2.03095$
and $\Im(\tilde{r}_{NT_{+}}(s_{0},\mu))$ diverging as predicted with
leading order asymptotics more and more closely resembling $-\frac{1}{2\pi}(s_{0}-\frac{1}{2})\frac{1}{\mu}$.

As noted this in turn provides a strong justification for ascribing
value $0$ to the 2d Cesaro limit of $\ln\left(\frac{z}{\tilde{z}}\right)$
when $\mu=0$ (as done in {[}1{]}) in order to retain continuity of
analytic continuation across $\mu=0$; and presumably similarly for
$\mu=-1$ and $\mu=-2$.

\section{Appendices}

\subsection{Tools for performing numerical calculations used in this paper}

The numerical results presented in this paper are generated in R-scripts
({}``R code Gamma.R'' and {}``R code Zeta.R'') and an XL2007 spreadsheet
with VBA code ({}``RootIdentitiesZeta\_mu\_-1To0\_TestsB.xlsm'').
The supporting files for the zeta R-script are two .txt files ({}``First10000Primes.txt
containing the first $10,000$ primes and {}``zero\_all2m.txt''
containing the first $2,000,000$ $NT_{+}$ roots of $\zeta$ as obtained
from the zip file contained at {[}2{]}). These are also used to populate
the required areas of the XL spreadsheet in order for its functions
to work correctly. The spreadsheet has been saved without this data
embedded for reasons of space, which is also why the file for the
$2,000,000$ $NT_{+}$ roots is left for the reader to unzip and rename
as {}``zero\_all2m.txt''. For both the R-scripts and the XL spreadsheet
the code is public and completely transparent. It is clear from this
code how to run the applications (e.g how to have the functions correctly
read the input roots and primes) and for the XL spreadsheet there
is additional brief guidance given within the spreadsheet itself.

Both applications (the R-scripts and the XL/VBA) can also easily be
adapted in order to conduct further testing as desired.

\end{document}